\documentclass[11pt,reqno]{amsart}

\usepackage{hyperref}
\usepackage{graphicx}
\usepackage{color}
\usepackage{amsfonts,amssymb}
\usepackage{bbm} 
\usepackage{psfrag} 

\usepackage{mathrsfs}
\usepackage{a4wide}

\newtheorem{neu}{}[section]
\newtheorem{Cor}[neu]{Corollary}
\newtheorem*{Cor*}{Corollary}
\newtheorem{Thm}[neu]{Theorem}
\newtheorem*{Thm*}{Theorem}
\newtheorem{Prop}[neu]{Proposition}
\newtheorem*{Prop*}{Proposition}
\theoremstyle{definition}

\newtheorem*{Rmk*}{Remark}
\newtheorem{Rmk}[neu]{Remark}
\newtheorem{Ex}[neu]{Example}
\newtheorem*{Ex*}{Example}

\newtheorem*{Qu*}{Question}

\newtheorem{Def}[neu]{Definition}

\newcommand{\R}{\mathbb{R}}
\newcommand{\C}{\mathbb{C}}

\newcommand{\pf}{\longrightarrow}

\newcommand{\hpf}{\hookrightarrow}

\newcommand{\im}{\mathrm{im\,}}

\newcommand{\om}{\omega}

\newcommand{\A}{\mathcal{A}}

\renewcommand{\P}{\mathcal{P}}

\newcommand{\F}{\mathcal{F}}

\newcommand{\D}{\mathbb{D}}

\newcommand{\E}{\mathcal{E}}
\newcommand{\K}{\mathcal{K}}

\renewcommand{\H}{\mathrm{H}}

\newcommand{\Ham}{\mathrm{Ham}}

\newcommand{\CF}{\mathrm{CF}}

\newcommand{\HF}{\mathrm{HF}}

\newcommand{\RFH}{\mathrm{RFH}}

\newcommand{\SH}{\mathrm{SH}}

\renewcommand{\A}{\mathcal{A}}
\newcommand{\p}{\partial}

\newcommand{\beq}{\begin{equation}}
\newcommand{\beqn}{\begin{equation}\nonumber}
\newcommand{\eeq}{\end{equation}}

\newcommand{\bea}{\begin{equation}\begin{aligned}}
\newcommand{\bean}{\begin{equation}\begin{aligned}\nonumber}
\newcommand{\eea}{\end{aligned}\end{equation}}

\newcommand{\HH}{{\widehat{H}}}

\numberwithin{equation}{section}

\definecolor{Urs}{rgb}{0,.7,0}
\definecolor{Peter}{rgb}{0,0,1}
\definecolor{red}{rgb}{1,0,0}

\begin{document}
\title[Non-displaceable contact embeddings and leaf-wise intersections]{Non-displaceable contact embeddings and infinitely many leaf-wise intersections}
\author{Peter Albers}
\author{Mark McLean}
\address{
    Peter Albers\\
    Department of Mathematics\\
    ETH Z\"urich, Switzerland}
\email{peter.albers@math.ethz.ch}
\address{
    Mark McLean\\
    Department of Mathematics\\
    ETH Z\"urich, Switzerland}
\email{mark.mclean@fim.math.ethz.ch}
\keywords{Symplectic homology, leaf-wise intersections, Lefschetz fibrations}
\subjclass[2000]{53D40, 37J10, 53D35}
\begin{abstract}
We construct using Lefschetz fibrations a large family of contact manifolds with the following properties:  Any bounding contact embedding into an exact symplectic manifold satisfying a mild topological assumption is non-displaceable and generically has infinitely many leaf-wise intersection points. Moreover, any Stein filling has infinite dimensional symplectic homology. 
\end{abstract}
\maketitle

\section{Main Results}\label{sec:main_result}

To state the main results we begin with the following construction of a large family of contact manifolds.
We denote by $\D^2(\delta):=\{z\mid|z|<\delta\}$ and by $\D^2$ the closed unit disk. Let $\pi:\widetilde{E}\pf\D^2$ be a Lefschetz fibration of dimension greater than $2$ with at least one critical point and fibers which are Liouville domains, see Definitions \ref{def:Liouville_domain} and \ref{def:convex_lef_fibr}. Without loss of generality we assume that $0\in\D^2$ is a regular value of $\pi$. We prove in Proposition \ref{prop:E_is_a_convex_lef_fibr_over_annulus} that $E:=\widetilde{E}\setminus\pi^{-1}(\D^2(\delta))$ can be made into a convex Lefschetz fibration over the annulus $A:=\D^2\setminus\D^2(\delta)$. After appropriate smoothing of the codimension two corners we obtain a contact manifold $(\Sigma,\xi)$. Let $(M,\om=d\lambda)$ be an exact symplectic manifold which is convex at infinity, that is, outside a compact set $M$ is symplectomorphic to the positive part of the symplectization of a compact contact manifold. 

\begin{Thm}\label{thm:main}
With the above notation let $\iota:\Sigma\hpf M$ be a contact embedding, that is, $\ker(\iota^*\lambda)=\xi$. Moreover, we assume that $\iota(\Sigma)$ bounds some compact region $V$ and that the canonical map $\H^1(V)\twoheadrightarrow\H^1(\Sigma)$ is surjective. For a generic embedding $\iota$ and a generic compactly supported Hamiltonian diffeomorphism $\phi\in\Ham_c(M,\om)$ there exist infinitely many leaf-wise intersections (with respect to $\phi$ and $\Sigma$).
\end{Thm}

We recall that $x\in\iota(\Sigma)$ is a leaf-wise intersection with respect to $\phi\in\Ham_c(M,d\lambda)$ if $\phi(x)\in L_x$ where $L_x\subset \iota(\Sigma)$ is the leaf of the Reeb flow of $\lambda|_\Sigma$ through $x\in\iota(\Sigma)$.

\begin{Cor} Under the same assumptions as in Theorem \ref{thm:main}, $\Sigma$ admits no embedding as a displaceable, bounding contact hypersurface into any exact symplectic manifold.
\end{Cor}

\begin{proof}
Suppose for a contradiction that $\iota(\Sigma)$ is displaced by a Hamiltonian $\phi$. Since this is an open condition we may assume that  $\iota$ and $\phi$ are generic. This contradicts Theorem \ref{thm:main} since leaf-wise intersection points are in particular intersection points $\iota(\Sigma)\cap\phi\big(\iota(\Sigma)\big)$.
\end{proof}

\begin{Thm}\label{thm:main2}
Under the same assumptions as in Theorem \ref{thm:main}, any strong symplectic filling $V$ of $\Sigma$ has infinite dimensional symplectic homology.
\end{Thm}

\begin{Rmk}
It follows from the proof of Theorem \ref{thm:main} that the above construction can be generalized to convex Lefschetz fibrations over a surface with boundary such that the monodromy map around one boundary component is the identity. 
\end{Rmk}

\begin{Rmk}
The homological condition $\H^1(V)\twoheadrightarrow\H^1(\Sigma)$ is too strong. It is only used in the proof of Proposition \ref{prop:infinite_SH} where we need that a certain closed closed one form extends. Moreover, this homological condition is automatically satisfied if the bounded region is a Stein domain.
\end{Rmk}

\begin{Ex}
Let $\widetilde{E}$ be an affine variety which is not a product. A closed embedding $\widetilde{E}\hpf\C^N$ induces a Lefschetz fibration $\widetilde{E}\pf\C$ with fiber being $\widetilde{E}\cap H$ where $H\subset\C^N$ is a generic hypersurface. Then $\Sigma$ is the boundary of the natural Liouville domain whose completion is $\widetilde{E}\setminus H$.

Another source of examples come from closed integral symplectic manifolds by removing two transverse Donaldson hypersurfaces (induced from the same ample line bundle).
\end{Ex}

\section*{Acknowledgement}
\noindent This article was written when both authors were at ETH Z\"urich. We thank ETH Z\"urich and in particular the symplectic working group for the stimulating working atmosphere.

This article has benefited enormously from numerous discussions with Urs Frauenfelder to whom we are indebted.

The first author is partially supported by NSF grant DMS-0805085.

\subsection*{Idea of the proof} To prove Theorem \ref{thm:main}, we first prove Theorem \ref{thm:main2}. This is done by showing that there is an infinite dimensional subgroup in symplectic homology generated by Reeb orbits of $\Sigma$ contained in a particular set of homotopy classes, see Theorem \ref{thm:nice_at_infinity_has_big_SH}. These Reeb orbits are concentrated near the fiber over zero in $\widetilde{E}$. We then use work by \cite{Cieliebak_Frauenfelder_Oancea_Rabinowitz_Floer_homology_and_symplectic_homology} and \cite{Albers_Frauenfelder_Leafwise_intersections_and_RFH,Albers_Frauenfelder_Leafwise_Intersections_Are_Generically_Morse} to conclude Theorem \ref{thm:main}.

\subsection*{Organizations of the article} In Section \ref{sec:main_definitions} we recall the basic definitions which enter in our construction. In Section \ref{sec:nice_at_infty} we introduce the notion of being nice at infinity and state the crucial Theorem \ref{thm:nice_at_infinity_has_big_SH}. This theorem is proved in Section \ref{sec:proof_of_crucial_thm}. In Section \ref{sec:rephrasing} we show that the above examples are in fact nice at infinity. Finally, in Section \ref{sec:final_proofs} we prove Theorems \ref{thm:main} and \ref{thm:main2}.

\section{Main Definitions}\label{sec:main_definitions}

In this section we recall some definitions and properties. We rely mostly on \cite{McLean_Lefschetz_fibrations_and_symplectic_homology}.

\begin{Def} \label{def:Liouville_domain}
A compact exact symplectic manifold $(M,d\lambda)$ with boundary is called a \textit{Liouville domain} if $(\p M,\lambda)$ is a contact manifold and if the vector field $Z$ defined by the equation $\iota_Zd\lambda=\lambda$ is transverse to $\p M$ and pointing outward. $Z$ is called the \textit{Liouville vector field}. 
\end{Def}
%


\begin{Def}
The \textit{completion} of a Liouville domain $(M,d\lambda)$ is the symplectic manifold
\beq
 \widehat{M}:=M\cup_{\p M}(\p M\times [1,\infty))
\eeq
where we extend the symplectic form by $d(r\lambda|_{\p M})$, $r\in [1,\infty)$, over the cylindrical end.
\end{Def}

\begin{Def}
A \textit{Liouville deformation} from $(M,d\lambda_0)$ to $(\widetilde{M},d\widetilde{\lambda})$ is a smooth family $(M,d\lambda_t)$, $t\in[0,1]$, such that $(M,d\lambda_1)$ is exact symplectomorphic to $(\widetilde{M},d\widetilde{\lambda})$. 
\end{Def}

Next we sketch the definition of symplectic homology $\SH_*(M,\lambda)$ of a Liouville domain $(M,d\lambda)$, see \cite{Oancea_survey_SH} for more details. We assume that the contact form $\lambda|_{\p M}$ is non-degenerate. For $a\in\R$ let $\HH^a:\widehat{M}\to\R$ be a function such that on the cylindrical end it satisfies $\HH^a(x,r)=ar$. For generic $a$ and for a small time-dependent perturbation $\HH^a_t$, Hamiltonian Floer homology $\HF_*(\HH^a_t)$ is well defined. The underlying complex $\CF_*(\HH^a_t,J)$ is generated by periodic orbits $\P(\HH^a_t)$ of $\HH^a_t$ and the differential is defined by counting rigid perturbed $J$-holomorphic cylinders. 

When $a<b$ there is a natural map $\HF_*(\HH^a_t)\to\HF_*(\HH^b_t)$. Symplectic homology is by definition the direct limit
\beq
\SH_*(M,\lambda):=\varinjlim_a\HF_*(\HH^a_t)\;.
\eeq
It is independent of all choices and  is an invariant of $\widehat{M}$ up to exact symplectomorphism, see \cite[Section 7b]{McLean_Lefschetz_fibrations_and_symplectic_homology}.

\begin{Def}\label{def:convex_lef_fibr}
A smooth map $\pi:E\pf S$ is a \textit{convex Lefschetz fibration} with Liouville form $\Theta$ and fiber $F$ if the following holds: (colon added)
\begin{enumerate}
\item $E$ is a manifold with boundary and codimension 2 corners such that the boundary can be decomposed as $\p E=\p^hE\cup\p^vE$.
\item $S$ is a compact 2-dimensional manifold with boundary.
\item The map $\pi$ satisfies the following conditions: 
\begin{enumerate}
\item $\pi$ has finitely many critical points of a certain local model, see \cite[Definition 2.12]{McLean_Lefschetz_fibrations_and_symplectic_homology} for details.
\item Outside the critical set $\pi$ is a submersion with typical fiber $F$, where $F$ is a smooth manifold with boundary.
\item There exists a open set $N^h\subset E$ and open neighborhood $N^{\p F}\subset F$ of $\p F$ such that
\beq
 E|_{N^h}:=N^h \cong S\times N^{\p F}
\eeq
as a fiber bundle with respect to $\pi$  and the projection map to $S$.
\end{enumerate} 
\item The vertical boundary is given by $\p^vE=E|_{\p S}:=\pi^{-1}(\p S)$ and $\p^v E\pf\p S$ is a fiber bundle. 
\item There exists a 1-form $\Theta$ on $E$ such that 
\begin{enumerate}
\item $d\Theta$ is a symplectic form,
\item The Liouville vector field $Z$ is transverse to $\p E$ and pointing outwards.
\item $(F,\theta_F:=\Theta|_F)$ is a Liouville domain.
\item On $N^h$ we have
\beq
\Theta|_{N^h}=\pi^*(sd\vartheta)+{pr}_2^*\theta_F
\eeq
where $pr_2$ is the projection $pr_2:S\times N^{\p F}\pf N^{\p F}$.
\end{enumerate}
\end{enumerate}
\end{Def}

In \cite[Section 2.2]{McLean_Lefschetz_fibrations_and_symplectic_homology} it has been shown that a convex Lefschetz fibration $\pi:E\to S$ with fiber $F$ admits a completion $\widehat{\pi}:\widehat{E}\to\widehat{S}$ with fiber $\widehat{F}$.  Moreover, \cite[Section 2.4]{McLean_Lefschetz_fibrations_and_symplectic_homology} has proven that in the definition of symplectic homology $\SH_*(\widehat{E})$ one can use functions of the form $\widehat{H}^a=\widehat{\pi}^*H_{\widehat{S}}^a+\widehat{pr}_2^*H^a_{\widehat{F}}$ on the cylindrical end of $\widehat{E}$ where $\widehat{pr}_2:\widehat{S}\times \big(\p F\times[1,\infty)\big)\to\p F\times[1,\infty)$ is the natural extension of $pr_2$ from the previous definition.

The symplectic form $d\Theta$ induces a connection on $E$ by taking the $d\Theta$-orthogonal plane field to the vertical tangent spaces. The parallel transport maps of this connection preserve the symplectic form. The monodromy map associated to a loop $\gamma:S^1\to S$ avoiding the critical values of $\pi$ in the base is the symplectomorphism $F_{\gamma(0)}\to F_{\gamma(0)}$ induced by parallel transport around $\gamma$.

\section{Exact symplectic manifolds nice at infinity}\label{sec:nice_at_infty}

\begin{Def} 
We define the annulus $\A:=[-1,1]\times S^1$ with coordinates $(s,\vartheta)$.
\end{Def}

\begin{Def}\label{def:being_compatible_with_fibration}
Let  $\pi:E\pf \A$ be a convex Lefschetz fibration with Liouville form $\Theta$ and fiber $F$. We say $\pi$ is \textit{trivial over an end} if the following is satisfied: ( colon added)
\begin{enumerate}
\item On $U_+:=(-\epsilon,1]\times S^1$  we require 
\beq
E|_{U_+}\cong U_+\times F
\eeq
as a fiber bundle with respect to $\pi$. 
\item For $\Theta$ we require that on $U_+$ we have
\beq
\Theta|_{U_+}=\pi^*(sd\vartheta)+pr_2^*\theta_F
\eeq
where $pr_2$ is by abuse of notation the projection $pr_2:U_+\times F\pf F$.
\end{enumerate}
\end{Def}

\begin{Def}\label{def:nice_at_infinity}
Let $(\E,d\Lambda)$ be an exact symplectic manifold with boundary with codimension 2 corners and let $\K$ be a compact subset of $\E$. Then we call the triple $(\E,\Lambda,\K)$ \textit{nice at infinity} with fiber $F$ if the following holds: There exists a convex Lefschetz fibration $\pi:E\pf\A$ with Liouville form $\Theta$ and fiber $F$ which is trivial over an end. Moreover, there is a compact subset $K\subset E$ with the following properties: ( colon added)
\begin{enumerate}
\item $K$ does not intersect the sets $N^h$ and $\pi^{-1}(U_+)$.
\item $\E\setminus\K$ is exact symplectomorphic to $E\setminus K$ via a symplectomorphism $\Psi$.
\item The 1-form $\Psi^*\left((\pi^*d\vartheta)|_{E\setminus K}\right)$ extends to a closed 1-form $\beta$ on $\E$.
\end{enumerate}
\end{Def}

\begin{Rmk}
After smoothing the corners of $\E$ it becomes a Liouville domain. In particular, we can associate symplectic homology $\SH_*(\E,\Lambda)$ to it.
\end{Rmk}

\begin{Thm}\label{thm:nice_at_infinity_has_big_SH}
Let $(\E,\Lambda,\K)$ be nice at infinity with fiber $F$. If  $\SH_*(F,\theta_F)\neq0$  then
\beq
\dim\SH_*(\E,\Lambda)=\infty\;.
\eeq
\end{Thm}
\section{Proof of Theorem \ref{thm:nice_at_infinity_has_big_SH}}\label{sec:proof_of_crucial_thm}

We start with some preliminary considerations. 

\begin{Rmk}
Let $(F,d\theta_F)$ be a compact symplectic manifold with convex boundary. Then $F$ admits a completion $\widehat{F}:=F\cup_{\p F}\big(\p F\times[1,\infty)\big)$ with one form $\theta_{\widehat{F}}:=r\cdot\theta_F|_{\p F}$ on $F\times[1,\infty)$ where $r$ is the radial coordinate.

Let $(E,d\Theta)$ be as in Definition \ref{def:being_compatible_with_fibration} with fiber $F$. Then in \cite[Section 2.2]{McLean_Lefschetz_fibrations_and_symplectic_homology} the completion $\widehat{\pi}:\widehat{E}\pf\widehat{\A}=\R\times S^1$ of $E$ is defined and has the following properties.
\begin{enumerate}
\item On the region $(-\epsilon,\infty)\times S^1$ we have 
\beq\label{eqn:prop1_of_E}
\Big(\widehat{E}|_{(-\epsilon,\infty)\times S^1},\Theta\Big)\cong\Big((-\epsilon,\infty)\times S^1\times\widehat{F},\widehat{\pi}^*(sd\vartheta)+pr_2^*(\theta_{\widehat{F}})\Big)
\eeq
\item On the region $(-\infty, -2)\times S^1$ we have 
\beq
\Big(\widehat{E}|_{(-\infty,-2)\times S^1},\Theta\Big)\cong\big((-\infty,-2)\times M(\phi),\theta_\phi\big)
\eeq
where $(M(\phi),\theta_\phi)$ is the mapping torus of the monodromy map $\phi$ around the loop $\{-1\}\times S^1\subset\A$, see last paragraph in Section \ref{sec:main_definitions} for the definition of monodromy.
Moreover, it contains the trivial bundle 
\beq
\big((-\infty,-2)\times S^1\times(\widehat{F}\setminus F),\widehat{\pi}^*(sd\vartheta)+pr_2^*(\theta_{\widehat{F}})\big)
\eeq
 as a subbundle.
\item On the region $(-2,-\epsilon)\times S^1$ the bundle $\widehat{E}$ contains the trivial subbundle 
\beq
\big((-2,-\epsilon)\times S^1\times(\widehat{F}\setminus F),\widehat{\pi}^*(sd\vartheta)+pr_2^*(\theta_{\widehat{F}})\big)
\eeq
which extends the previous trivial bundle in the obvious manner.
\end{enumerate}
\end{Rmk}

\begin{figure}[!h]
\psfrag{F}{$F$}
\psfrag{as}{$as$}
\psfrag{-eps}{$-\epsilon$}
\psfrag{-eps/4}{$-\frac{\epsilon}{4}$}
\psfrag{-1}{$-2$}
\psfrag{-1+eps}{$-2+\epsilon$}
\psfrag{1/2s}{$\frac12 s$}
\psfrag{Anu}{$\A$}
\psfrag{E}{$E$}
\psfrag{K}{$K$}
\psfrag{X}{$Z$}
\psfrag{Nh}{$N^h$}
\psfrag{pi}{$\pi$}
\psfrag{Gamma}{$\Gamma$}
\psfrag{EU-}{$E|_{U_-}$}
\psfrag{Up x F}{$E|_{U_+}\cong U_+\times F$}
\psfrag{Nhcong}{$N^h\cong\A\times N^{\F}$}
 \includegraphics[scale=.95]{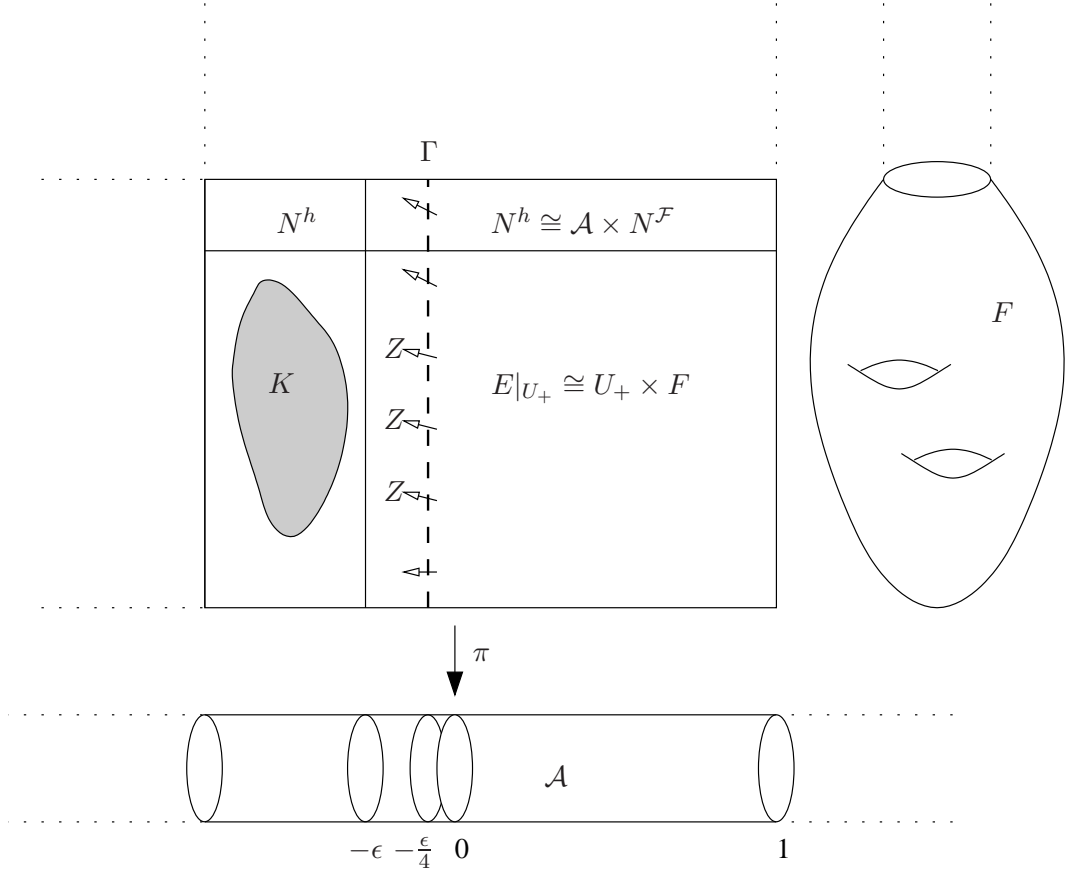}
 \caption{\label{fig:bundle_E} The completion $\widehat{E}$. The objects $Z$, $\Gamma$, etc.~will be introduced later.}
\end{figure}

\begin{Def}
Let $(\E,d\Lambda)$ be as in Definition \ref{def:nice_at_infinity}.  Then we define the completion $\widehat{\E}$ as
\beq
\widehat{\E}:=(\widehat{E}\setminus K)\cup\K
\eeq
\end{Def}

\begin{figure}[!h]
\psfrag{H_F_hat}{$H^a_{\widehat{F}}$}
\psfrag{slope a}{slope $=a$}
\psfrag{F}{$F$}
\psfrag{dF}{$\p F$}
\psfrag{dFc}{$\p F\times[1,\infty)$}
 \includegraphics[scale=0.75]{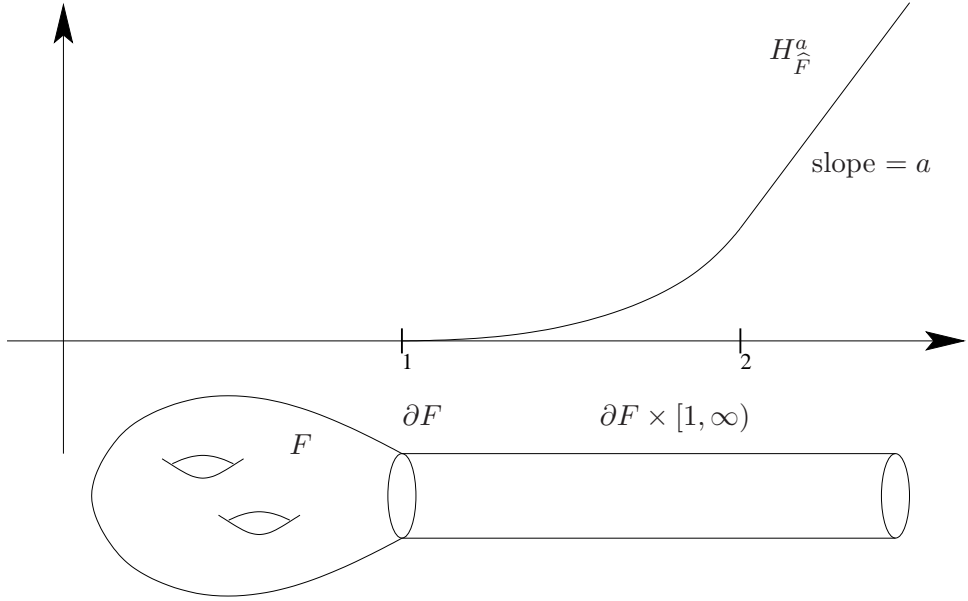}
 \caption{\label{fig:H_F_hat} The function $H^a_{\widehat{F}}$.}
\end{figure}

On the completion $\widehat{F}$ for $a\geq1$ we define the  function
\beq
H^a_{\widehat{F}}(x):=\begin{cases}
0 & x\in F\\
f_a(r) & x=(y,r)\in\p F\times[1,2]\\
ar & x=(y,r)\in\p F\times[2,\infty]\\
\end{cases}
\eeq
where $f_a$ is a smooth function with $f_a'\geq0$, $f_a''\geq0$ making $H^a_{\widehat{F}}$ into a smooth function, see Figure \ref{fig:H_F_hat}. 

Similarly, we define a function $H^a_0:S^1\times\R\pf\R$ by
\beq
H^a_0(s,t):=\begin{cases}
-as & s\in(-\infty,-2]\\
f_a(-s) & s\in [-2,-1]\\ 
0 & s\in[-1,-\epsilon]\\
g_a(s) & s\in[-\epsilon,1]\\
as & s\in[1,\infty]
\end{cases}
\eeq
where $g_a$ is a smooth function with $g_a'\geq0$, $g_a''\geq0$ making $H^a_0$ into a smooth function. Moreover, we require
\beq\label{eqn:property_of_g_a}
g_a'(s)=\tfrac12\quad s\in[-\epsilon/2,0]
\eeq
We choose $g'_a(s)=\frac12$ on $[-\epsilon/2,0]$ to ensure that $H^a_0$ has no periodic orbits in this region.

\begin{figure}[!h]
\psfrag{H_0}{$H^a_0$}
\psfrag{-as}{$-as$}
\psfrag{as}{$as$}
\psfrag{-eps}{$-\epsilon$}
\psfrag{-eps/2}{$-\frac{\epsilon}{2}$}
\psfrag{1/2s}{$\frac12 s$}
 \includegraphics[scale=.75]{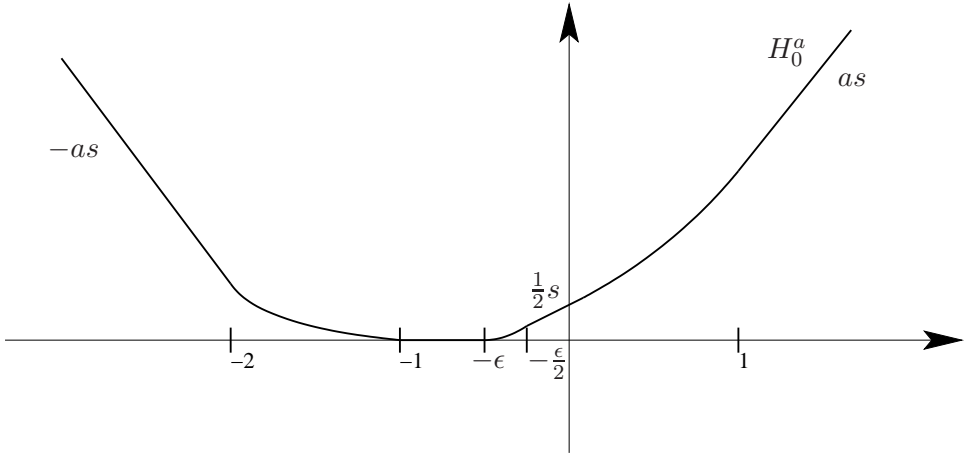}
 \caption{\label{fig:H_0} The function $H^a_0$.}
\end{figure}

\begin{Def}\label{def:Hamiltonian_adapted_to_completion}
We call a function $\widehat{H}:\widehat{E}\pf\R$ \textit{adapted to the completion} if it satisfies
\begin{enumerate}
\item On the region $\widehat{E}|_{(-\epsilon,\infty)\times S^1}$ we have 
\beq\label{eqn:direct_sum_Ham}
\widehat{H}=\widehat{\pi}^*H_0^a+pr_2^*H^a_{\widehat{F}}.
\eeq
\item On the region $\widehat{E}|_{(-\infty, -2)\times S^1}$ we have 
\beq
\widehat{H}=\widehat{\pi}^*H_0^a+pr_2^*H^a_{\widehat{F}}
\eeq
which makes sense since $H^a_{\widehat{F}}=0$ in the region of $\widehat{E}$ where the fibration is non-trivial.
\item On the region $\widehat{E}|_{(-2,-\epsilon)\times S^1}$  we have
\beq
\widehat{H}=pr_2^*H^a_{\widehat{F}}
\eeq
with the same understanding as above.
\item Everywhere else we choose $\widehat{H}$ such that it has only constant periodic orbits which are non-degenerate.
\end{enumerate}
This definition naturally extends to functions $\widehat{H}:\widehat{\E}\pf\R$ since the regions mentioned (1) -- (3) are disjoint from $\K$ (cf. to Definition \ref{def:nice_at_infinity} for the notation).

The real number $a\geq1$ is called the asymptotic slope of $\widehat{H}$.
\end{Def}

\begin{proof}[Proof of Theorem \ref{thm:nice_at_infinity_has_big_SH}]
Let $(\E,d\Lambda)$ be nice at infinity with fiber $F$ and completion $\widehat{\pi}:\widehat{\E}\pf\widehat{\A}=\R\times S^1$. After a small perturbation of $\p F$ we may assume that the contact manifold $(\p F,\theta_F|_{\p F})$ has only non-degenerate Reeb orbits. In particular, it has discrete period spectrum, i.e.~the set $\{T\in\R\mid T\text{ is a the period of a Reeb orbit}\}$ is a discrete subset of $\R$. 

We choose a function $\widehat{H}^a:\widehat{E}\pf\R$ which is adapted to the completion. On the region mentioned in (1) -- (3) in Definition \ref{def:Hamiltonian_adapted_to_completion} the Hamiltonian function $\widehat{H}^a$ is degenerate. For an open and dense set of asymptotic slopes the only degeneracy of $\widehat{H}^a$ comes from the fact that $\widehat{H}^a$ is autonomous since $(\partial F,\theta_F|_{\p F})$ is non-degenerate. Then after a small time-dependent perturbation localized near the periodic orbits we obtain a non-degenerate Hamiltonian function $\widehat{H}^a_t$. This can be done so that on the region $\widehat{E}|_{(-\epsilon,\infty)\times S^1}$ the Hamiltonian function $\widehat{H}^a_t$ is still a sum of two Hamiltonian functions as in equation \eqref{eqn:direct_sum_Ham} but now with time-dependent  summands:
\beq
\widehat{H}^a_t(\cdot)=\widehat{\pi}^*H_0^a(t,\cdot)+pr_2^*H^a_{\widehat{F}}(t,\cdot).
\eeq
Moreover, for a sufficiently small perturbation the perturbed orbits remain in the same homotopy class as the unperturbed orbit.
 
Thus, after the perturbation the non-constant periodic orbits are contained in $\widehat{\E}\setminus \K$. Those orbits contained in the region $\widehat{E}|_{(-\epsilon/4,\infty)\times S^1}$ project via $\widehat{\pi}$ to loops in $\R\times S^1$ with negative winding (measured with respect to $d\vartheta$). The orbits contained in the region $\widehat{E}|_{(-\infty, -2)\times S^1}$ project to loops with positive winding.
There are no orbits in the region $\widehat{E}|_{[-\epsilon/2,0]\times S^1}$ by equation \ref{eqn:property_of_g_a}. Finally, the  orbits contained in the region $\widehat{E}|_{(-2,-\epsilon/4)\times S^1}$ project to contractible loops in $\R\times S^1$. 

We divide the set $\P(\widehat{H}^a_t)$ of periodic orbits of $\widehat{H}^a_t$ as follows
\beq
\P^\text{I}:=\left\{x\in\P(\widehat{H}^a_t)\left|\;\int_{S^1} x^*\beta >0\right.\right\}
\eeq
where $\beta$ is the 1-form from Definition \ref{def:nice_at_infinity}. We set $\P^\text{II}:=\P(\widehat{H}^a_t)\setminus\P^\text{I}$. We define $\CF_\text{I}^a$ to be the sub-complex of the symplectic homology complex $\CF(\widehat{H}^a_t)$ generated by $\P^\text{I}$. Analogously, we define $\CF_\text{II}^a$. In fact, Stokes theorem implies that the condition $\int_{S^1} x^*\beta >0$ is preserved under the Floer differential, thus $\CF(\widehat{H}^a_t)$ splits as a direct sum 
\beq\label{eqn:SH_splits}
\CF(\widehat{H}^a_t)=\CF_\text{I}^a\oplus\CF_\text{II}^a\;.
\eeq
We denote by $\HF(\widehat{H}^a_t)=\HF_\text{I}^a\oplus\HF_\text{II}^a$ the corresponding homology groups. To finish the proof of Theorem \ref{thm:nice_at_infinity_has_big_SH} we will show that $\dim \HF_\text{I}^a=\infty$.

For this we show that 
\beq\label{eqn:equaliyt_PI}
\P^\text{I}=\{x\in\P(\widehat{H}^a_t)\mid x(t)\in\widehat{E}|_{(-\epsilon/4,\infty)\times S^1}\;\forall t\in S^1\}\;.
\eeq
Indeed, since $\int_{S^1} x^*\beta>0$ the orbit $x$ is non-constant, thus $x$ is contained in $\widehat{\E}\setminus\K$. Thus, we compute using notation from Definition \ref{def:nice_at_infinity}
\bea
0<\int_{S^1} x^*\beta&=\int_{S^1} x^*\Psi^*\left((\pi^*d\vartheta)|_{E\setminus K}\right)\\
&=\int_{S^1} (\pi\circ\Psi\circ x)^*d\vartheta=\text{winding of the projection of }x
\eea
Thus, equality \eqref{eqn:equaliyt_PI} follows.

Next we prove that any Floer cylinder $u:\R\times S^1\pf\widehat{E}$ connecting two periodic orbits in $\P^\text{I}$ is entirely contained in $\widehat{E}|_{(-\epsilon/4,\infty)\times S^1}$. This follows from the maximum principle \cite[Lemma 7.2]{Abouzaid_Seidel_Viterbo_functoriality}. We consider the Liouville vector field $Z$ defined by the equation $d\Lambda(Z,\cdot)=\Lambda$. We claim that $Z$ is transversal to the hypersurface $\Gamma:=\{-\frac{\epsilon}{4}\}\times S^1\times F$ and points into the region $\widehat{\E}\setminus\left(\widehat{E}|_{(-\epsilon/4,\infty)\times S^1}\right)$, see Figure \ref{fig:bundle_E}. Assuming this claim for the moment it follows immediately from \cite[Lemma 7.2]{Abouzaid_Seidel_Viterbo_functoriality} that 
\beq
\im(u)\cap\left(\widehat{\E}\setminus\left(\widehat{E}|_{(-\epsilon/4,\infty)\times S^1}\right)\right)\subset \Gamma\;.
\eeq 
Thus, it remains to prove that $Z$ has the claimed properties. Since $(\widehat{\E}|_{(-\epsilon,0)\times S^1},\Lambda)$ is (exact symplectomorphic to) a product and by equation \eqref{eqn:prop1_of_E} the Liouville vector field $Z$ projects to the Liouville vector field of the 1-form $sd\vartheta$ on $(-\epsilon,0)\times S^1$. Since $s<0$ the Liouville vector field $Z$ has the required properties.

Since any Floer cylinder connecting two periodic orbits in $\P^\text{I}$ is entirely contained in $\widehat{E}|_{(-\epsilon/4,\infty)\times S^1}$ we can compute $\HF_\text{I}^a$ inside 
\beq
\left(\widehat{E}|_{(-\epsilon/4,\infty)\times S^1}, \Theta\right)\cong\left(\big((-\epsilon/4,\infty)\times S^1\big)\times \widehat{F},\widehat{\pi}^*(sd\vartheta)+pr_2^*(\theta_{\widehat{F}})\right)
\eeq
with respect to the Hamiltonian function $\widehat{H}^a_t(\cdot)=\widehat{\pi}^*H_0^a(t,\cdot)+pr_2^*H^a_{\widehat{F}}(t,\cdot)$. Thus, 
\beq\label{eqn:computation_SH_I}
\HF_\text{I}^a=\HF(K_0^a(t,\cdot))\otimes\HF(H^a_{\widehat{F}}(t,\cdot))
\eeq
on the Liouville domain $\big((\R\times S^1)\times\widehat{F},sd\vartheta+\theta_{\widehat{F}}\big)$ where $K_0^a:S^1\times(\R\times S^1)\pf\R$ is defined as 
\beq
K_0^a(t,s,\vartheta)=
\begin{cases}
H_0^a(t,s,\vartheta) & s\geq-\frac{\epsilon}{4}\\[1ex]
\frac12 s& \text{else}
\end{cases}
\eeq
Now, taking the direct limit over $a\to\infty$ we obtain from $\HF_\text{I}^a$ and equation \eqref{eqn:computation_SH_I}
\beq
\SH_\text{I}=G\otimes \HF(F,\theta_F)
\eeq
where $G$ is the part homology of loop space of $S^1$ consisting of loops with positive winding number, see \cite{Viterbo_partII, Salamon_Weber_Floer_homology_and_heat_flow, Abbo_Schwarz_On_the_Floer_homology_of_cotangent_bundles}. In particular, if $\SH(F,\theta_F)\neq0$
\beq
\dim \SH_\text{I}=\infty\;.
\eeq
This concludes the proof.
\end{proof}

\section{Rephrasing the construction}\label{sec:rephrasing}

Let $\pi:\widetilde{E}\pf\D^2$ be a convex Lefschetz fibration and $\D(\delta):=\{z\in\C\mid|z|<\delta\}$. We assume without loss of generality that $0\in\D^2$ is a regular value. 

\begin{Prop}\label{prop:E_is_a_convex_lef_fibr_over_annulus}
Let $E:=\widetilde{E}\setminus\pi^{-1}(\D^2(\delta))$, then the map $\pi:E\pf\A\cong(\D^2\setminus\D^2(\delta))$ can be given the structure of a convex Lefschetz fibration which is trivial over an end.
\end{Prop}
  
\begin{proof}
We assume that (after a Liouville deformation) the Lefschetz fibration $\widetilde{E}$ is trivial over $\D^2(2\delta)$, that is
\beq
\big(\widetilde{E}|_{\D^2(2\delta)},\widetilde{\Theta}|_{\D^2(2\delta)}\big)\cong \big(\D^2(2\delta)\times F, \pi^*(\tfrac12 r^2d\varphi)+{pr}_2^*\theta_F\big)
\eeq
where $(r,\varphi)$ are polar coordinates on $\D^2$. We choose a smooth function $\rho:[0,1]\to\R$ such that
\beq
\rho(r)=\begin{cases}
-1 & r=\delta\\
0 & r=\tfrac32\delta\\
1 & r=2\delta\\
\end{cases}
\eeq
and such that $\rho'>0$. Then for $K\gg0$ sufficiently large the 1-form $\Theta$ on $E$ defined by
\beq
\Theta:=\widetilde{\Theta}+K\pi^*\big(\rho(r)d\varphi\big)
\eeq
is a Liouville form, see Definition \ref{def:convex_lef_fibr}. Then for a suitable $\delta<\delta'<\frac32\delta$ the region $\widetilde{E}\setminus\pi^{-1}(\D^2(\delta'))$ is a convex Lefschetz fibration with a trivial end. Shifting the coordinate $r$ proves the proposition.
\end{proof}

\section{Proof of Theorems \ref{thm:main} and \ref{thm:main2} }\label{sec:final_proofs}

We use the notation from Theorem \ref{thm:main}.

\begin{Prop}\label{prop:infinite_SH}
The Liouville domain $(V,d\lambda)$ has infinite dimensional symplectic homology
\beq
\dim\SH_*(V,\lambda)=\infty\;.
\eeq
\end{Prop}

\begin{proof}
Recall that $(\Sigma,\alpha)$ is a contact type hypersurface in the exact symplectic manifold $(E,d\alpha)$ obtained by smoothing the  codimension two corners of $\p E$. Furthermore, $(\Sigma,\iota^*\lambda)$ is the contact manifold induced by the embedding into $M$ with the same contact structure: $\xi=\ker\alpha=\ker\lambda$. Thus, (after rescaling $\alpha$) there exists an exact symplectic cobordism with negative end $(\Sigma,\iota^*\lambda)$ and positive end $(\p E,\alpha)$. We attach $V$ to the negative end of this cobordism creating an exact symplectic manifold $(\E,d\Lambda)$ with codimension 2 corners which is nice at infinity with fiber $F$ according to Definition \ref{def:nice_at_infinity} where condition (3) in Definition \ref{def:nice_at_infinity} is satisfied due to the assumption  $\H^1(V)\twoheadrightarrow\H^1(\Sigma)$ in Theorem \ref{thm:main}. This is illustrated in Figure \ref{fig:E_tilde}.

\begin{figure}[!h]
\psfrag{Eh}{$\widehat{E}$}
\psfrag{sigma}{$(\Sigma,\alpha$)}
\psfrag{sigma2}{$(\Sigma,\iota^*\lambda$)}
\psfrag{cob}{$\mathrm{cobordism}$}
\psfrag{p E}{$\p E$}
\psfrag{V}{$V$}
 \includegraphics[scale=.75]{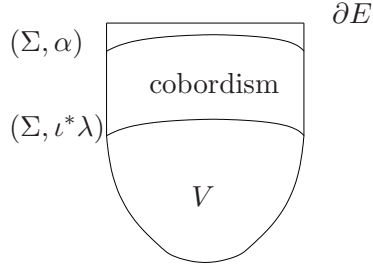}
 \caption{\label{fig:E_tilde} The exact symplectic manifold $(\E,d\Lambda)$.}
\end{figure}

By assumption the Lefschetz fibration $\widetilde{E}$ has at least one critical point, and thus each fiber contains an exact Lagrangian sphere, the vanishing cycle. It follows immediately from \cite[Theorem 4.3]{Viterbo_partI} that $\SH_*(F)\neq0$. Thus, by Theorem \ref{thm:nice_at_infinity_has_big_SH} 
\beq
\dim\SH_*(\E,\Lambda)=\infty\;.
\eeq
Since $(\widehat{\E},\Lambda)$ is exact symplectomorphic to $(\widehat{V},\lambda)$, they have isomorphic symplectic homology groups. 
\end{proof}

\begin{proof}[Proof of Theorem \ref{thm:main}]
Theorems 1.2 and 1.5 in \cite{Cieliebak_Frauenfelder_Oancea_Rabinowitz_Floer_homology_and_symplectic_homology} give us a long exact sequence 
\beq
\cdots\to\SH^{-*}(V)\stackrel{a}{\to}\SH_*(V)\to\RFH_*(V,\Sigma)\to\cdots\;.
\eeq
Proposition 1.3 in \cite{Cieliebak_Frauenfelder_Oancea_Rabinowitz_Floer_homology_and_symplectic_homology} then asserts that the map $a$ factors through some finite dimensional vector space. Thus,  $\dim\SH_*(V)=\infty$ implies that $\dim \RFH_*(V,\Sigma)=\infty$. Finally, Proposition 3.1 in \cite{Cieliebak_Frauenfelder_Oancea_Rabinowitz_Floer_homology_and_symplectic_homology} states that
\beq
\RFH_*(V,\Sigma)\cong\RFH_*(M,\Sigma)\;.
\eeq
It follows from \cite[Proposition 2.4, Theorem 2.13]{Albers_Frauenfelder_Leafwise_intersections_and_RFH} and \cite[Theorem 2.5]{Albers_Frauenfelder_Leafwise_Intersections_Are_Generically_Morse} that for generic $\iota$ and $\phi$ the number of leaf-wise intersections is at least as big as $\dim\RFH_*(V,\Sigma)=\infty$.
\end{proof}

\begin{proof}[Proof of Theorem \ref{thm:main2}]
Let $V$ be a strong filling of $\Sigma$ such that the canonical map $\H^1(V)\twoheadrightarrow\H^1(\Sigma)$ is surjective. Then $\Sigma$ has a contact embedding into the completion $\widehat{V}$ of $V$ satisfying all the assumptions of Theorem \ref{thm:main}. Thus, Proposition \ref{prop:infinite_SH} applies.
\end{proof}

%
%
\bibliographystyle{amsalpha}
\bibliography{../../../Bibtex/bibtex_paper_list}
\end{document}